\documentclass[a4,psfig,english,11pt,epsf,portrait]{article}
\usepackage{amsmath,amsfonts,latexsym,amscd,amssymb,theorem}
\usepackage{epsfig}
\usepackage{amsmath}
\def\R{\hbox{\bf R}}
\def\Z{\hbox{\bf Z}}

\def\N{\hbox{\bf N}}

\def\g{\gamma}
\def\b{\beta}

\def\O{\Omega}

\def\vphi{\varphi}

\def\<{\langle}
\def\>{\rangle}
\newcommand{\ba}{\begin{eqnarray}}
\newcommand{\ea}{\end{eqnarray}}

\newtheorem{thm}{Theorem}[section]

\newtheorem{theorem}[thm]{Theorem}
\newtheorem{definition}[thm]{Definition}
\newtheorem{lemma}[thm]{Lemma}
\newtheorem{proposition}[thm]{Proposition}

\newtheorem{rem}[thm]{Remark}

\newcommand{\eps}{\epsilon}
\numberwithin{equation}{section}

\renewcommand{\R}{{\mathbb R}}
\renewcommand{\Z}{{\mathbb Z}}
\renewcommand{\N}{{\mathbb N}}

\renewcommand{\g}{\gamma}

\begin{document}

\title{\bf A critical parabolic Sobolev embedding via Littlewood-Paley decomposition}
\author{
\normalsize\textsc{ H. Ibrahim \footnote{CEREMADE, Universit\'{e}
Paris-Dauphine, Place De Lattre de Tassigny, 75775 Paris Cedex 16,
France. E-mail: ibrahim@cermics.enpc.fr
\newline
\indent $\,\,{}^{1}$LaMA-Liban, Lebanese University, P.O. Box 826
Tripoli, Lebanon
} $^{ ,\,1}$}}
\vspace{20pt}

\maketitle

\centerline{\small{\bf{Abstract}}} \noindent{\small{In this paper,
we show a parabolic version of the Ogawa type inequality in Sobolev
spaces. Our inequality provides an estimate of the $L^{\infty}$ norm
of a function in terms of its parabolic $BMO$ norm, with the aid of
the square root of the logarithmic dependency of a higher order
Sobolev norm. The proof is mainly based on the Littlewood-Paley
decomposition and a characterization of parabolic $BMO$ spaces.}}

\hfill\break
 \noindent{\small{\bf{AMS subject classifications:}}} {\small{42B35,
     54C35, 42B25, 39B05.}}\hfill\break
  \noindent{\small{\bf{Key words:}}} {\small{Littlewood-Paley
      decomposition, logarithmic
      Sobolev inequalities, parabolic $BMO$ spaces, Lizorkin-Triebel
      spaces, Besov spaces.}}\hfill\break


\section{Introduction and main results}\label{sec1}
In order to study the long-time existence of a certain class of singular
parabolic problems, Ibrahim and Monneau \cite{IM09} made use of a
parabolic logarithmic Sobolev inequality. They proved that for
$f\in W^{2m,m}_{2}(\R^{n+1})$, $m,n \in \N^{*}$ and $2m>\frac{n+2}{2}$,
the following estimate takes place (with $\log^{+} x =
\max(\log x, 0)$):
\begin{equation}\label{eq1:IM}
\|f\|_{L^{\infty}(\R^{n+1})} \leq C (1+
  \|f\|_{BMO^{a}(\R^{n+1})}(1+
    \log^{+}\|f\|_{W^{2m,m}_{2}(\R^{n+1})} )),
\end{equation}
for some constant $C=C(m,n)>0$. Here $BMO^{a}$ stands for the
anisotropic Bounded Mean Oscillation space with the parabolic
anisotropy $a=(1,\ldots,1,2)\in \R^{n+1}$ (see
Definition~\ref{BMOa}), while $W^{2m,m}_{2}$ stands for the
parabolic Sobolev space (see Definition~\ref{para_sob}). The above
estimate, after also being proved on a bounded domain
\begin{equation}\label{omega_T}
\Omega_{T} = (0,1)^{n} \times (0,T) \subseteq \R^{n+1},
\end{equation}
was successfully applied in order to obtain some \textit{a priori}
bounds on the gradient of the solution of particular parabolic
equations leading eventually to the long-time existence (see
\cite[Proposition 3.7]{IM09} or \cite[Theorem 1.3]{IJM08}). The
bounded version of (\ref{eq1:IM}) (see \cite[Theorem 1.2]{IM09})
reads: if $f\in W^{2m,m}_{2}(\Omega_T)$ with $2m>\frac{n+2}{2}$,
then:
\begin{equation}\label{Leb:1}
\|f\|_{L^{\infty}(\Omega_T)}\leq C (1+
\|f\|_{\overline{BMO}^{a}(\Omega_T)} (1+ \log^{+}
\|f\|_{W^{2m,m}_{2}(\Omega_T)})),
\end{equation}
where $C=C(m,n,T)>0$ is a positive constant, and
\begin{equation}\label{bar_BMO}
\|f\|_{\overline{BMO}^{a}(\Omega_{T})} = \|f\|_{BMO^{a}(\Omega_{T})}
+ \|f\|_{L^{1}(\Omega_{T})}.
\end{equation}
Indeed, the fact that inequality (\ref{eq1:IM}) does not hold on
$\O_{T}$ with a positive constant $C^*=C^{*}(m,n,T)$ can be easily
understood by applying this inequality to the function
$f=(C^{*}+\eps) \in W^{2m, m}_{2}(\O_{T})$ with $\eps>0$. In this
case $\|f\|_{L^{\infty}(\O_{T})} = C^{*} + \eps$,
$\|f\|_{BMO^{a}(\O_{T})} = 0$, and hence a contradiction. However,
working on $\R^{n+1}$, the same function $f$ could not be used since
$f\in\!\!\!\!\!/\,W^{2m, m}_{2}(\R^{n+1})$. Let us indicate that
both inequalities (\ref{eq1:IM}) and (\ref{Leb:1}) still hold for
vector-valued functions $f=(f_{1},\ldots,f_{n},f_{n+1})\in
(W^{2m,m}_{2}(\R^{n+1}))^{n+1}$ with $2m>\frac{n+2}{2}$ and the
natural change in norm.

The elliptic version of (\ref{eq1:IM}) was showed by Kozono and
Taniuchi in \cite{KT00}. Indeed, they have showed that for $f\in
W_{p}^{s}(\R^{n})$, $1<p<\infty$, the following estimate holds:
\begin{equation}\label{eq2:KT}
\|f\|_{L^{\infty}(\R^{n})} \leq C (1 + \|f\|_{BMO(\R^{n})}(1 +
\log^{+}\|f\|_{W^{s}_{p}(\R^{n})})),\quad sp>n,
\end{equation}
for some $C=C(n,p,s)>0$. Here $BMO$ is the usual elliptic/isotropic
bounded mean oscillation space (defined via Euclidean balls). The
main advantage of (\ref{eq2:KT}) is that it was successfully applied
in order to extend the blow-up criterion of solutions to the Euler
equations originally given by Beale, Kato and Majda in \cite{BKM84}.
This blow-up criterion was then refined by Kozono, Ogawa and
Taniuchi \cite{KOT02}, and by Ogawa \cite{Og03}, showing weaker
regularity criterion that was even relaxed by Planchon
\cite{Planchon03}, Danchin \cite{Danchin05}, and Cannone, Chen and
Miao \cite{CCM07}.

The proof of inequality (\ref{eq1:IM}) is based on the analysis in
anisotropic Lizorkin-Triebel, Besov, Sobolev and $BMO^{a}$ spaces.
This is made via Littlewood-Paley decomposition and various Sobolev
embeddings. In fact, some of the technical arguments were inspired
by Ogawa \cite{Og03} in his proof of the sharp version of
(\ref{eq2:KT}) that reads:  if $g\in L^{2}(\R^{n})$ and $f:=\nabla g
\in W_{q}^{1}(\R^{n})\cap L^{2}(\R^{n})$ for $n<q$, then there
exists a constant $C=C(q)>0$ such that:
\begin{equation}\label{eq3:O}
\|f\|_{L^{\infty}(\R^{n})} \leq C(q)\left(1 + \|f\|_{BMO(\R^{n})}
\left(\log^{+}(\|f\|_{ W_{q}^{1}(\R^{n})} +
  \|g\|_{L^{\infty}(\R^{n})})\right)^{1/2}\right).
\end{equation}
It is worth mentioning that the original type of the logarithmic
Sobolev inequalities (\ref{eq2:KT}) and (\ref{eq3:O}) was found in
Br\'{e}zis and Gallou\"{e}t \cite{BG80}, and Br\'{e}zis  and Wainger
\cite{BW80}. The Br\'{e}zis-Gallou\"{e}t-Wainger inequality states
that the $L^{\infty}$ norm of a function can be estimated by the
$W_{p}^{n/p}$ norm with the partial aid of the $W_{r}^{s}$ norm with
$s>n/r$ and $1\leq r\leq \infty$. Precisely,
\begin{equation}\label{cpMSW}
\|f\|_{L^{\infty}(\R^{n})} \leq C\Big((1 + \log (1 +
\|f\|_{W_{r}^{s}(\R^{n})}))\Big)^{\frac{p-1}{p}}
\end{equation}
holds for all $f\in W_{p}^{n/p}(\R^{n})\cap W_{r}^{s}(\R^{n})$ with
the normalization $\|f\|_{W_{p}^{n/p}(\R^{n})} = 1$. Originally,
Br\'{e}zis and Gallou\"{e}t \cite{BG80} obtained (\ref{cpMSW}) for
the case $n=p=r=s=2$, where they applied their inequality in order
to prove global existence of solutions to the nonlinear
Schr\"{o}dinger equation. Later on, Br\'{e}zis and Wainger
\cite{BW80} obtained (\ref{cpMSW}) for the general case, and
remarked that the power $\frac{p-1}{p}$ in (\ref{cpMSW}) is optimal
in the sense that one can not replace it by any smaller power.
However, it seems that little is known about the sharp constant in
(\ref{cpMSW}).

Coming back to inequalities (\ref{eq1:IM}), (\ref{eq2:KT}) and
(\ref{eq3:O}), the natural question that arises is the following:
why does the inequality (\ref{eq1:IM}) seems to be the parabolic
extension of (\ref{eq2:KT}) although the proof is inspired (as
mentioned above) from that of (\ref{eq3:O}) given by Ogawa
\cite{Og03}? The answer to this question is partially contained in
\cite[Remark 2.14]{IM09} where the authors pointed out that the
well-known relation between elliptic/isotropic Lizorkin-Triebel and
$BMO$ spaces (see \cite[Proposition 2.3]{Og03}) will not be used in
the proof of (\ref{eq1:IM}) even though it seems to be valid
(without giving a proof) in the parabolic/anisotropic framework. The
relation is the following:
\begin{equation}\label{eq4:equiv}
\dot{F}^{0,a}_{\infty,2} \simeq BMO^{a},
\end{equation}
where $\dot{F}^{0,a}_{\infty,2}$ is the homogeneous parabolic
Lizorkin-Triebel space (see Definition~\ref{LT_defi}).

In this paper, we show a parabolic version of the logarithmic
Sobolev inequality (\ref{eq3:O}) basically using the equivalence
(\ref{eq4:equiv}) that is shown to be true (see Lemma~\ref{lem1}).
This answers the question raised above. Our study takes place on the
whole space $\R^{n+1}$ and on the bounded domain $\O_T$. A
comparison (in some special cases) of our inequality with
(\ref{eq1:IM}) is also discussed.

Before stating our main results, we define some terminology. A
generic element in $\R^{n+1}$ will be denoted by $z=(x,t)\in
\R^{n+1}$ where $x=(x_{1},\ldots,x_{n})\in \R^{n}$ is the spatial
variable, and $t\in \R$ is the time variable. For a given function
$g$, the notation $\partial_{i}g$ stands for the partial derivative
with respect to the spatial variable: $\partial_{i}g =
\partial_{x_{i}} g := \frac{\partial g}{\partial x_{i}}$, $i=1,...,
n$. In this case $\partial_{n+1}g = \partial_{t} g:= \frac{\partial
g}{\partial
  t}$. We also denote $\partial^{s}_{x} g$, $s\in \N$, any derivative with
respect to $x$ of order $s$. Moreover, we denote the space-time
gradient by $\nabla g := (\partial_{1}g, \ldots,
\partial_{n}g, \partial_{n+1}g)$. Finally, we denote $\|f\|_{X}:=
\max (\|f_{1}\|_{X},\ldots,\|f_{n}\|_{X}, \|f_{n+1}\|_{X})$ for any
vector-valued function $f=(f_{1},\ldots,f_{n}, f_{n+1})\in X^{n+1}$
where $X$ is any Banach space. Throughout this paper and for the
sake of simplicity, we will drop the superscript $n+1$ from
$X^{n+1}$. Following the above notations, our first theorem reads:

\begin{theorem}\textit{(Parabolic Ogawa inequality on $\R^{n+1}$)}.\label{theorem} Let $m, n\in \N^*$
with
$2m>\frac{n+2}{2}$. Then there exists a constant $C=C(m,n)>0$ such
that for any function $g\in L^{2}(\R^{n+1})$ with
$f=(f_{1},\ldots,f_{n},f_{n+1})=\nabla g
\in~W^{2m,m}_{2}(\R^{n+1})$, we have:
\begin{equation}\label{I_inq}
\|f\|_{L^{\infty}(\R^{n+1})} \leq C\left(1 +
\|f\|_{BMO^{a}(\R^{n+1})}\left(\log^{+}(\|f\|_{
W_{2}^{2m,m}(\R^{n+1})} +
  \|g\|_{L^{\infty}(\R^{n+1})})\right)^{1/2}\right).
\end{equation}
\end{theorem}

\begin{rem}\label{rem1}
All the terms appearing in (\ref{I_inq}) make sense since for
$2m>\frac{n+2}{2}$, there exists some $\g=\g(m,n)>0$ such that:
$$W^{2m,m}_{2} \hookrightarrow C^{\g,\g/2} \hookrightarrow
L^{\infty}\hookrightarrow BMO^{a},\quad 0<\g<1,$$ where
$C^{\g,\g/2}$ is the usual parabolic H\"{o}lder space. Moreover, it
is easy to see that $g$ is continuous and bounded.
\end{rem}

\begin{rem}\label{rem_s}
By taking $m,n\in \N^{*}$, $2m>\frac{n+2}{2}$, the same inequality
(\ref{I_inq}) holds for $g\in L^{\infty} (\R^{n+1})$ and $ f =
\partial_{i} g \in W^{2m,m}_{2}(\R^{n+1})$ for some fixed
$i=1,\ldots,n+1$. This can be considered as the scalar-valued
version of the vector-valued version (\ref{I_inq}).
\end{rem}

\begin{rem}\label{rem2}
Inequalities (\ref{eq1:IM}) and (\ref{I_inq}) have the same order of the
higher regular term. As a consequence, inequality (\ref{I_inq}) can also
be applied in order to establish the long-time existence of solutions of
the parabolic problems studied in \cite{IJM08, IM09}.
\end{rem}
Our next theorem concerns a similar type inequality of
(\ref{I_inq}), but with functions $g$ and $f$ defined over $\O_T$
(given by (\ref{omega_T})). Before stating this result, we first
remark that in the case of functions $f=\nabla g$ defined on a
bounded domain, we formally have (by Poincar\'e inequality):
$$\|g\|_{L^{\infty}}\leq C \|f\|_{L^{\infty}},$$
where $C>0$ is a constant depending on the measure of the domain.
Moreover, since
$$\|f\|_{L^{\infty}} \leq C_{1} \|f\|_{C^{\g,\g/2}}\leq C_{2}
\|f\|_{W^{2m,m}_{2}}\quad\mbox{with}\quad C_{1}, C_{2}>0,$$ the
above two estimates imply that the term $\|g\|_{L^{\infty}}$ should
be dropped from inequality (\ref{I_inq}) when dealing with functions
defined over bounded domains. Indeed, we have:
\begin{theorem}\textit{(Parabolic Ogawa inequality on a bounded
domain)}.\label{theorem-bdd} Let $f\in W^{2m,m}_{2}(\O_T)$ with
$2m>\frac{n+2}{2}$. Then there exists a constant $C=C(m,n,T)>0$ such
that:
\begin{equation}\label{leb:2}
\|f\|_{L^{\infty}(\O_{T})} \leq C \left(1 +
\|f\|_{\overline{BMO}^{a}(\O_T)}\left(\log^{+}\|f\|_{
W_{2}^{2m,m}(\O_T)}\right)^{1/2}\right),
\end{equation}
where the norm $\|\cdot\|_{\overline{BMO}^{a}(\O_{T})}$ is given by
(\ref{bar_BMO}).
\end{theorem}
\begin{rem}\label{rem_sh}
Inequality (\ref{leb:2}) is sharper than (\ref{Leb:1}) by the simple
observation that $x^{1/2} \leq 1 + x$ for $x\geq 0$. In other words,
inequality (\ref{leb:2}) implies (\ref{Leb:1}) with the same
positive constant $C = C(m,n)$.
\end{rem}
In the same spirit of Remark~\ref{rem_sh}, our last theorem gives a
comparison between inequality (\ref{eq1:IM}) and (\ref{I_inq}) for a
certain class of functions $g$, and for particular space dimensions.
\begin{theorem}\textit{(Comparison between parabolic logarithmic
inequalities)}. \label{theorem3} Let $n=1,2,3$ and $m\in \N^{*}$
satisfying $2m>\frac{n+2}{2}$. There exists a constant $C=C(m,n)>0$
such that for the class of functions $g\in L^{2}(\R^{n+1})$ with
$\|g\|_{L^{2}(\R^{n+1})} \leq 1$, and $f=\nabla g \in
W^{2m,m}_{2}(\R^{n+1})$, we have:
\begin{equation}\label{sec5:eq1}
\big(\log^{+}(\|f\|_{W^{2m,m}_{2}(\R^{n+1})} +
\|g\|_{L^{\infty}(\R^{n+1})})\big)^{1/2}\leq C(1 +
\log^{+}\|f\|_{W^{2m,m}_{2}(\R^{n+1})}),
\end{equation}
and hence inequality (\ref{I_inq}) implies (\ref{eq1:IM}) for
possibly a different positive constant $C$.
\end{theorem}
\subsection{Organization of the paper}
This paper is organized as follows. In Section \ref{sec2}, we
present some definitions and the main tools used in our analysis.
This includes parabolic Littlewood-Paley decomposition and various
Sobolev embeddings. Section~\ref{sec3} is devoted to the proof of
Theorem~\ref{theorem} (estimate on the entire space $\R^{n+1}$)
using mainly the equivalence (\ref{eq4:equiv}) that we also show in
Lemma~\ref{lem1}. In Section~\ref{sec4}, we give the proof of
Theorem~\ref{theorem-bdd} (estimate on the bounded domain $\O_T$).
Finally, in Section~\ref{sec5}, we give the proof of
Theorem~\ref{theorem3}.


\section{Preliminaries and basic tools}\label{sec2}
In this section, we define the fundamental function spaces used in
this paper. We also recall some important embeddings.
\subsection{Parabolic $BMO^{a}$ and Sobolev spaces}
Each coordinate $x_{i}$, $i=1, ..., n$ is given the weight $1$, while the
time coordinate $t$ is given the weight $2$. The vector
$a=(a_{1},\ldots, a_{n}, a_{n+1}) = (1,\ldots, 1, 2)\in \R^{n+1}$
is called the $(n+1)$-dimensional parabolic anisotropy. For
this given $a$, the action of $\mu\in [0,\infty)$ on $z=(x,t)$ is given
by $\mu^{a}z = (\mu x_{1},\ldots, \mu x_{n}, \mu^{2} t)$. For
$\mu>0$ and $s\in \R$ we set $\mu^{s a} z = (\mu^{s})^{a} z$. In
particular, $\mu^{-a} z = (\mu^{-1})^{a} z$ and $2^{-j a} z =
(2^{-j})^{a} z$, $j\in \Z$. For $z\in \R^{n+1}$, $z\neq 0$, let
$|z|_{a}$ be the unique positive number $\mu$ such that:
\begin{equation*}
\frac{x^{2}_{1}}{\mu^{2}}+\cdots+\frac{x^{2}_{n}}{\mu^{2}}+\frac{t^{2}}{\mu^{4}}
= 1
\end{equation*}
and let $|z|_{a} = 0$ for $z=0$. The map $|\cdot|_{a}$ is called the
parabolic distance function which is $C^{\infty}$ (see for instance
\cite{Y86}). In the case where $a=(1,\ldots, 1)\in \R^{n+1}$, we get
the usual Euclidean distance
$\|z\|=(x^{2}_{1}+\cdots+x^{2}_{n}+t^{2})^{1/2}$. Denoting
$\mathcal{O}\subseteq \R^{n+1}$, any open subset of $\R^{n+1}$, we
are ready to give the definition of the first two parabolic spaces
used in our analysis.
\begin{definition}\textit{(Parabolic bounded mean oscillation spaces)}.\label{BMOa}
A function $f\in L^{1}_{loc}(\mathcal{O})$ (defined up to an
additive constant) is said to be of parabolic bounded mean
oscillation, $f\in BMO^{a}(\mathcal{O})$, if we have:
\begin{equation}\label{s3mwdo3}
\|f\|_{BMO^{a}(\mathcal{O})} = \sup_{\mathcal{Q}\subseteq
\mathcal{O}} \inf_{c\in \R} \left(
  \frac{1}{|\mathcal{Q}|}\int_{\mathcal{Q}} |f-c| \right)<+\infty,
\end{equation}
where $\mathcal{Q}$ denotes (for $z_{0}\in \mathcal{O}$ and $r>0$)
an arbitrary parabolic cube:
$$\mathcal{Q} = \mathcal{Q}_{r}(z_{0}) = \{z\in \R^{n+1};\;
|z-z_{0}|_{a}<r\}.$$
\end{definition}
\begin{definition}\textit{(Parabolic Sobolev spaces)}.\label{para_sob}
Let $m\in \N$. We define the parabolic Sobolev space
$W^{2m,m}_{2}(\mathcal{O})$ as follows:
$$W^{2m,m}_{2}(\mathcal{O}) = \{f\in L^{2}(\mathcal{O});\;
\partial^{r}_{t}\partial_{x}^{s} f\in L^{2}(\mathcal{O}), \forall r,s\in
\N \mbox{ such that }2r+s\leq 2m\},$$ with
$\|f\|_{W^{2m,m}_{2}(\mathcal{O})} =
\sum^{2m}_{j=0}\sum_{2r+s=j}\|\partial^{r}_{t}\partial_{x}^{s}
f\|_{L^{2}(\mathcal{O})}$.
\end{definition}

\subsection{Parabolic Lizorkin-Triebel and Besov spaces}
Along with the above parabolic distance $|\cdot|_{a}$, the
Littlewood-Paley decomposition is now recalled (for more details, we
refer to \cite{Grafakos}). Let $\theta\in C_{0}^{\infty}(\R^{n+1})$
be any cut-off function satisfying:
\begin{equation}\label{old_neta}
\theta(z) = \left\{
\begin{aligned}
& 1 \quad &\mbox{if}& \quad |z|_{a}\leq 1\\
& 0 \quad &\mbox{if}& \quad |z|_{a}\geq 2.
\end{aligned}
\right.
\end{equation}
Let $\psi(z) = \theta(z) - \theta(2^{a} z)$. We now construct a
smooth (compactly supported) parabolic dyadic partition of unity
$(\psi_{j})_{j\in \Z}$ by letting
\begin{equation}\label{dyadic}
\psi_{j}(z) = \psi(2^{-ja} z),\quad j\in \Z,
\end{equation}
satisfying
$$\sum_{j\in \Z} \psi_{j}(z) = 1\quad \mbox{for}\quad z\neq 0.$$
Define $\vphi_{j}$, $j\in \Z$, as the inverse Fourier transform of
$\psi_{j}$, i.e. $\hat{\vphi}_{j} = \psi_{j}$ where we let
\begin{equation}\label{y5tche}
\vphi:=\vphi_{0}.
\end{equation}
It is worth noticing that $\vphi_{j}$ satisfies:
\begin{equation}\label{lp:eq2}
\vphi_{j}(z) = 2^{(n+2) j} \vphi(2^{ja} z),\quad j\in \Z\quad
\mbox{and}\quad z\in \R^{n+1}.
\end{equation}
The above Littlewood-Paley decomposition asserts that any tempered
distribution $f\in \mathcal{S}'(\R^{n+1})$ can be decomposed as:
$$f = \sum_{j\in \Z} \vphi_{j} * f\quad \mbox{with the convergence in }
\mathcal{S}'/\mathcal{P} \mbox{ (modulo polynomials)}.$$ Here
$\mathcal{S}(\R^{n+1})$ is the usual Schwartz class of rapidly
decreasing functions and $\mathcal{S}'(\R^{n+1})$ is its
corresponding dual, represents the space of tempered distributions.
We now define parabolic Lizorkin-Triebel spaces.

\begin{definition}\textit{(Parabolic homogeneous Lizorkin-Triebel
spaces)}.\label{LT_defi} Given a smoothness parameter $s\in \R$, an
integrability exponent $1\leq p<\infty$, and a summability exponent
$1\leq q\leq \infty$. Let $\vphi_j$ be given by (\ref{lp:eq2}), we
define the parabolic homogeneous Lizorkin-Triebel space
$\dot{F}^{s,a}_{p,q}$ as the space of all functions $f\in
\mathcal{S}'(\R^{n+1})$ with finite quasi-norms
$$
\|f\|_{\dot{F}^{s,a}_{p,q}(\R^{n+1})} = \left\|\left( \sum_{j\in \Z}
2^{s q j} |\vphi_{j} * f|^{q} \right)^{1/q}
\right\|_{L^{p}(\R^{n+1})}<\infty,
$$
and the natural modification for $q=\infty$, i.e.
$$
\|f\|_{\dot{F}^{s,a}_{p,\infty}(\R^{n+1})} = \Big\|\sup_{j\in \Z}
2^{s j} |\vphi_{j} * f| \Big\|_{L^{p}(\R^{n+1})}.
$$
In the case $p=\infty$ and $s=0$, we define the parabolic
homogeneous Lizorkin-Triebel space $\dot{F}^{0,a}_{\infty,q}$ as the
space of all functions $f\in \mathcal{S}'(\R^{n+1})$ with finite
quasi-norms:
$$\|f\|_{\dot{F}^{0,a}_{\infty,q}} = \sup_{\mathcal{Q}\in P} \left(\frac{1}{|\mathcal{Q}|}
\int_{\mathcal{Q}}\sum_{j=-\mbox{scale}(\mathcal{Q})}^{\infty}
|\vphi_{j}*f|^{q} \right)^{1/q}<\infty,$$ where $P$ is the
collection of all dilated parabolic cubes $\mathcal{Q} = 2^{a j}
[(0,1)^{n+1} + k]$, with $\mbox{scale}(\mathcal{Q}) = j\in \Z$ and
$k\in \Z^{n+1}$.
\end{definition}
As a convention, for $s\in \R$, and $1\leq q<\infty$, we denote
\begin{equation}\label{conv1}
\|f_{+}\|_{\dot{F}^{s,a}_{\infty,q}(\R^{n+1})} = \Big\|( \sum_{j\geq
1} 2^{s q j} |\vphi_{j} * f|^{q} )^{1/q}
\Big\|_{L^{\infty}(\R^{n+1})}
\end{equation}
and
\begin{equation}\label{conv2}
\|f_{-}\|_{\dot{F}^{s,a}_{\infty,q}(\R^{n+1})} = \Big\|( \sum_{j\leq
-1} 2^{s q j} |\vphi_{j} * f|^{q} )^{1/q}
\Big\|_{L^{\infty}(\R^{n+1})}.
\end{equation}
The space $\dot{F}^{0,a}_{p,2}$ can be identified with the parabolic
Hardy space $H^{p,a}(\R^{n+1})$, $1\leq p<\infty$, having the
following square function characterization stated informally as:
\begin{equation}\label{Hardy}
H^{p,a}(\R^{n+1}) = \Big\{f\in \mathcal{S}'(\R^{n+1});\; (\sum_{j\in
    \Z} |\vphi_{j}*f|^{2})^{1/2}\in L^{p}\Big\}.
\end{equation}
This identification between the above two spaces is the following:
\begin{theorem}\textit{(Identification between $H^{p,a}$ and
    $\dot{F}^{0,a}_{p,2}$)}. \textit{(See Bownik
      \cite{Bow07}.)}\label{theorem_bow1}
For all $1\leq p<\infty$, we have $\dot{F}^{0,a}_{p,2}(\R^{n+1})
\simeq H^{p,a}(\R^{n+1})$.
\end{theorem}
Another useful space throughout our analysis is the parabolic
inhomogeneous Besov space. The main difference in defining this
space is the choice of the parabolic dyadic partition of unity that
is now altered. Indeed, we take $(\psi_{j})_{j\geq 0}$ satisfying:
\begin{equation}\label{inhomo_dyadic}
\psi_{j}:= \left\{
\begin{aligned}
& \psi_{j}\;&\mbox{ defined by (\ref{dyadic}) if }&\; j\geq
1\\
& \theta \;&\mbox{ defined by (\ref{old_neta}) if }&\; j=0.
\end{aligned}
\right.
\end{equation}
Again, it is clear that $\sum_{j\geq 0} \psi_{j}(z) = 1$, but now
for all $z\in \R^{n+1}$, and in exactly the same way as above, we
can rewrite the Littlewood-Paley decomposition with
\begin{equation}\label{ran_ah}
\hat{\vphi}_{j} = \psi_{j},\; j\geq 0,\quad \mbox{$\psi_{j}$ is
given by (\ref{inhomo_dyadic})}.
\end{equation}
We then arrive to the following definition:
\begin{definition}\textit{(Parabolic inhomogeneous Besov
spaces)}.\label{inhpB} Given a smoothness parameter $s\in \R$, an
integrability exponent $1\leq p\leq\infty$, and a summability
exponent $1\leq q\leq \infty$, we define the parabolic inhomogeneous
Besov space $B^{s,a}_{p,q}$ as the space of all functions $u\in
\mathcal{S}'(\R^{n+1})$ with finite quasi-norms
$$
\|u\|_{B^{s,a}_{p,q}} = \left( \sum_{j\geq 0} 2^{s q j} \|\vphi_{j}
* u\|^{q}_{L^{p}(\R^{n+1})}\right)^{1/q}<\infty,\quad \vphi_{j}\mbox{ is
given by (\ref{ran_ah})}
$$
and the natural modification for $q=\infty$, i.e.
\begin{equation}\label{inhpB_norm2}
\|u\|_{B^{s,a}_{p,\infty}} = \sup_{j\geq 0} 2^{s j} \|\vphi_{j}
* u\|_{L^{p}(\R^{n+1})},\quad \vphi_{j}\mbox{ is given by
(\ref{ran_ah})}.
\end{equation}
\end{definition}
For a detailed study of anisotropic Lizorkin-Triebel and Besov spaces,
we refer the reader to Triebel \cite{Tri1}.
\subsection{Embeddings of parabolic Besov and Sobolev spaces}
We present two embedding results from Johnsen and Sickel
\cite{Joh_Sic07}, and St\"{o}ckert \cite{Stockert82}.
\begin{theorem}\textit{(Embeddings of Besov spaces)}.\textit{(See
    Johnsen and Sickel \cite{Joh_Sic07}.)} \label{joh_sic}
Let $s,t\in\R$, $s>t$, and $1\leq p,r\leq \infty$ satisfy:
$s-\frac{n+2}{p} = t-\frac{n+2}{r}$. Then for any $1\leq q\leq \infty$
we have the following continuous embedding:
\begin{equation}\label{b_embed}
B^{s,a}_{p,q}(\R^{n+1}) \hookrightarrow B^{t,a}_{r,q}(\R^{n+1}).
\end{equation}
\end{theorem}
\begin{proposition}\textit{(Sobolev embeddings in Besov
    spaces)}.\textit{(See St\"{o}ckert \cite{Stockert82}.)}
  \label{stockert}
Let $m\in \N$, then we have:
\begin{equation}\label{s_embed}
W^{2m,m}_{2}(\R^{n+1}) \hookrightarrow B^{2m,a}_{2,\infty}(\R^{n+1}).
\end{equation}
\end{proposition}


\section{Proof of Theorem~\ref{theorem}}\label{sec3}
In this section we give the proof of Theorem~\ref{theorem}. We start
by showing the equivalence (\ref{eq4:equiv}) whose isotropic version
can be found in Triebel \cite{Tri}, and Frazier and Jawerth
\cite{FraJaw90}.
\begin{lemma}\textit{(Equivalence between $\dot{F}^{0,a}_{\infty,2}$ and
    $BMO^{a}$)}.\label{lem1}
We have $\dot{F}^{0,a}_{\infty,2}(\R^{n+1}) \simeq BMO^{a}(\R^{n+1})$. Precisely, there
exists a constant $C>0$ such that:
\begin{equation}\label{salta}
C^{-1}\|f\|_{\dot{F}^{0,a}_{\infty,2}}\leq \|f\|_{BMO^{a}} \leq C
\|f\|_{\dot{F}^{0,a}_{\infty,2}}.
\end{equation}
\end{lemma}
\textbf{Proof.} Using the result of Bownik \cite[Theorem
1.2]{Bow08}, we have the following duality argument (that can be
viewed as the anisotropic extension of the well-known isotropic
result of Triebel \cite{Tri}, and Frazier and Jawerth
\cite{FraJaw90}):
\begin{equation}\label{dual1}
\left(\dot{F}^{0,a}_{1,2}\right)' \simeq \dot{F}^{0,a}_{\infty,2},
\end{equation}
where $\left(\dot{F}^{0,a}_{1,2}\right)'$ stands for the dual space of
$\dot{F}^{0,a}_{1,2}$. Applying Theorem~\ref{theorem_bow1} with $p=1$ we
obtain:
\begin{equation}\label{dual2}
\dot{F}^{0,a}_{1,2} \simeq H^{1,a}.
\end{equation}
Using the description of the dual of anisotropic Hardy spaces
$H^{p,a}$ for $0<p\leq 1$ (see Bownik \cite[Theorem 8.3]{Bow03}), we
get:
\begin{equation}\label{Bow_dua}
\left(H^{p,a}\right)' = \mathcal{C}^{l}_{q,s}
\end{equation}
with the terms $p, l, q, s$ chosen such that:
\begin{equation}\label{lawmd3l7}
\left\{
\begin{aligned}
& l = \frac{1}{p} - 1,\\
& 1\leq \frac{q}{q-1} \leq \infty \quad \mbox{and} \quad p<
\frac{q}{q-1},\\
& s\in \N \quad \mbox{and} \quad s\geq \lfloor l \rfloor,\quad
\lfloor l\rfloor = \max \{n\in \Z;\, n\leq l\}.
\end{aligned}
\right.
\end{equation}
The function space $\mathcal{C}^{l}_{q,s}$, $l\geq 0$, $1\leq
q<\infty$ and $s\in \N$ (called the \textit{Campanato space}), is
the space of all $f\in L^{q}_{loc}(\R^{n+1})$ (defined up to
addition by $P\in \mathcal{P}_{s}$; the set of all polynomials in
$(n+1)$ variables of degree at most $s$) so that:
\begin{equation}\label{D7kl7}
\|f\|_{\mathcal{C}^{l}_{s,q}(\R^{n+1})} = \sup_{\mathcal{Q}\subseteq
\R^{n+1}} \inf_{P\in \mathcal{P}_{s}} |\mathcal{Q}|^{l} \left(
\frac{1}{|\mathcal{Q}|} \int_{\mathcal{Q}} |f - P|^{q}\right)^{1/q}
< \infty.
\end{equation}
Choosing $p=1$, $l=0$, $q=1$ and $s=0$, we can easily see that
conditions (\ref{lawmd3l7}) are all satisfied, and that (see
(\ref{D7kl7}) and (\ref{s3mwdo3})):
$$\mathcal{C}^{0}_{1,0} \simeq BMO^{a}.$$
This identification, together with (\ref{Bow_dua}), finally give:
\begin{equation}\label{dual3}
\left(H^{1,a}\right)' \simeq BMO^{a}.
\end{equation}
The proof then directly follows from (\ref{dual1}), (\ref{dual2}) and
(\ref{dual3}).$\hfill{\Box}$\\

A basic estimate is now shown in the following lemma.
\begin{lemma}\textit{(Logarithmic estimate in parabolic Lizorkin-Triebel
    spaces)}.\label{lem2}
Let $\g>0$ be a positive real number. Then, for $f\in
\dot{F}^{0,a}_{\infty,1}$ with
$\|f_{+}\|_{\dot{F}^{\g,a}_{\infty,2}(\R^{n+1})}$ and
$\|f_{-}\|_{\dot{F}^{-\g,a}_{\infty,2}(\R^{n+1})}$ are both finite,
there exists a constant $C=C(n,\g)>0$ such that:
\begin{equation}\label{tark}
\|f\|_{\dot{F}^{0,a}_{\infty,1}} \leq C \left(1 +
\|f\|_{\dot{F}^{0,a}_{\infty,2}} \left(\log^{+}
  (\|f_{+}\|_{\dot{F}^{\g,a}_{\infty,2}} + \|f_{-}\|_{\dot{F}^{-\g,a}_{\infty,2}}
  )\right)^{1/2}\right).
\end{equation}
\end{lemma}
\textbf{Proof.} We first indicate that the constant $C=C(n,\g)>0$
may vary from line to line in the proof which is divided into two steps.\\

\noindent {\bf Step 1} \textit{(First estimate on
$\|f\|_{\dot{F}^{0,a}_{\infty,1}}$)}. Let $N\in \N$, we compute
\begin{eqnarray*}
\|f\|_{\dot{F}^{0,a}_{\infty,1}} &\leq& \Big\|\sum_{j<-N}2^{\g j}2^{-\g j}|\vphi_{j} *
  f|\Big\|_{L^{\infty}} +\Big\|\sum_{|j|\leq N}|\vphi_{j} *
  f|\Big\|_{L^{\infty}} + \Big\|\sum_{j>N}2^{-\g j}2^{\g j}|\vphi_{j} *
  f|\Big\|_{L^{\infty}}\\
&\leq& C_{\g}2^{-\g N}\Big\|\Big(\sum_{j<-N} 2^{-2\g j}|\vphi_{j} *
  f|^{2}\Big)^{1/2}\Big\|_{L^{\infty}} + (2N+1)^{1/2}
  \Big\|\Big(\sum_{|j|\leq N}|\vphi_{j} *
  f|^{2}\Big)^{1/2}\Big\|_{L^{\infty}}\\
& &\hspace{0.3cm}+C_{\g}2^{-\g N}\Big\|\Big(\sum_{j>N} 2^{2\g j}|\vphi_{j} *
  f|^{2}\Big)^{1/2}\Big\|_{L^{\infty}}\\
&\leq&C_{\g}2^{-\g N}\|f_{-}\|_{\dot{F}^{-\g,a}_{\infty,2}} +
C(2N+1)^{1/2} \|f\|_{\dot{F}^{0,a}_{\infty,2}} + C_{\g}2^{-\g
  N}\|f_{+}\|_{\dot{F}^{\g,a}_{\infty,2}},
\end{eqnarray*}
with $C_{\g} = \left(\frac{1}{2^{2\g} -1} \right)^{1/2}$. As a
conclusion we may write
\begin{equation}\label{mt2ato}
\|f\|_{\dot{F}^{0,a}_{\infty,1}}\leq C \left((2N+1)^{1/2}
\|f\|_{\dot{F}^{0,a}_{\infty,2}} + 2^{-\g N}
(\|f_{+}\|_{\dot{F}^{\g,a}_{\infty,2}} +
\|f_{-}\|_{\dot{F}^{-\g,a}_{\infty,2}})\right).
\end{equation}
{\bf Step 2} \textit{(Optimization in $N$)}. We optimize (\ref{mt2ato}) in
$N$ by setting:
$$N=1 \quad \mbox{if} \quad \|f_{+}\|_{\dot{F}^{\g,a}_{\infty,2}} +
\|f_{-}\|_{\dot{F}^{-\g,a}_{\infty,2}}\leq 2^{\g}
\|f\|_{\dot{F}^{0,a}_{\infty,2}}.$$
Then it is easy to check (using (\ref{mt2ato})) that
\begin{equation}\label{nas1}
\|f\|_{\dot{F}^{0,a}_{\infty,1}}\leq C\|f\|_{\dot{F}^{0,a}_{\infty,2}}
\left(1+ \left(\log^{+} \frac{\|f_{+}\|_{\dot{F}^{\g,a}_{\infty,2}} +
\|f_{-}\|_{\dot{F}^{-\g,a}_{\infty,2}}}{\|f\|_{\dot{F}^{0,a}_{\infty,2}}}\right)^{1/2}\right).
\end{equation}
In the case where $\|f_{+}\|_{\dot{F}^{\g,a}_{\infty,2}} +
\|f_{-}\|_{\dot{F}^{-\g,a}_{\infty,2}}>2^{\g}
\|f\|_{\dot{F}^{0,a}_{\infty,2}}$, we take $1\leq \b <2^{\g}$ such that
$$N = N(\b)=\log_{2^{\g}}^{+}\left(\b\frac{\|f_{+}\|_{\dot{F}^{\g,a}_{\infty,2}} +
\|f_{-}\|_{\dot{F}^{-\g,a}_{\infty,2}}}{\|f\|_{\dot{F}^{0,a}_{\infty,2}}}
\right) - \frac{1}{2}\in \N.$$
In fact this is valid since the function $N(\b)$ varies continuously
from $N(1)$ to $N(2^{\g}) = 1 + N(1)$ on the interval $[1,2^{\g}]$. Using
(\ref{mt2ato}) with the above choice of $N$, we obtain:
\begin{eqnarray*}
\|f\|_{\dot{F}^{0,a}_{\infty,1}} &\leq& C \left[2^{1/2}\left(
    \log_{2^{\g}}^{+}\left(\b\frac{\|f_{+}\|_{\dot{F}^{\g,a}_{\infty,2}}
        +
\|f_{-}\|_{\dot{F}^{-\g,a}_{\infty,2}}}{\|f\|_{\dot{F}^{0,a}_{\infty,2}}}
\right)\right)^{1/2} \|f\|_{\dot{F}^{\g,a}_{\infty,2}} +
\frac{2^{\g/2}}{\b}\|f\|_{\dot{F}^{\g,a}_{\infty,2}}  \right]\\
&\leq&  C \left[\frac{2}{(\g\log2)^{1/2}}\left( \log^{+}\left(\frac{\|f_{+}\|_{\dot{F}^{\g,a}_{\infty,2}} +
\|f_{-}\|_{\dot{F}^{-\g,a}_{\infty,2}}}{\|f\|_{\dot{F}^{0,a}_{\infty,2}}}
\right)\right)^{1/2} \|f\|_{\dot{F}^{\g,a}_{\infty,2}} +
\frac{2^{\g/2}}{\b}\|f\|_{\dot{F}^{\g,a}_{\infty,2}}  \right],
\end{eqnarray*}
where for the second line we have used the fact that
$$\log^{+}\b < \log^{+}\frac{\|f_{+}\|_{\dot{F}^{\g,a}_{\infty,2}} +
\|f_{-}\|_{\dot{F}^{-\g,a}_{\infty,2}}}{\|f\|_{\dot{F}^{0,a}_{\infty,2}}}.$$
The above computations again imply (\ref{nas1}). By using the inequality:
\begin{equation}\label{nbbo2}
x\left(\log \left(e + \frac{y}{x} \right)\right)^{1/2} \leq \left\{
\begin{aligned}
& C(1+ x(\log(e+y))^{1/2})\quad &\mbox{for}&\quad 0<x\leq 1,\\
& Cx(\log (e+y))^{1/2}\quad &\mbox{for}&\quad x>1,
\end{aligned}
\right.
\end{equation}
in (\ref{nas1}), we directly arrive to our result. $\hfill{\Box}$\\

We now present the proof of our first main result.\\

\noindent \textbf{Proof of Theorem~\ref{theorem}.} First let us
mention that the constant $C=C(m,n)>0$ appearing in the following
proof may vary from line to line. We will show inequality
(\ref{I_inq}) in the scalar-valued version, i.e. by considering $f=
f_{i} = \partial_{i} g$ for some fixed $i=1,\ldots,n+1$. The
vector-valued version can then be easily deduced. The proof requires
estimating all the terms of inequality~(\ref{tark}). We start with
the obvious estimate (see (\ref{salta})):
\begin{equation}\label{result1}
\|f\|_{\dot{F}^{0,a}_{\infty,2}} \leq C \|f\|_{BMO^{a}}.
\end{equation}
The remaining terms will be estimated in the following three steps.\\

\noindent {\bf Step 1} \textit{(An upper bound on
$\|f_{+}\|_{\dot{F}^{\g,a}_{\infty,2}}$)}. Set $\eta = 2m - \frac{n+2}{2}
>0$. Choose $\g$ such that:
$$0<\g<\eta.$$
We compute (see (\ref{conv1})):
\begin{eqnarray}\label{BaHu}
\|f_{+}\|_{\dot{F}^{\g,a}_{\infty,2}} &=& \left\|\left( \sum_{j\geq 1} 2^{2 \g j}
|\vphi_{j} * f|^{2} \right)^{1/2}\right\|_{L^{\infty}}\nonumber\\
&\leq& C
  \sup_{j\geq 1} 2^{\eta j} \|\vphi_{j} * f\|_{L^{\infty}}
\end{eqnarray}
with $C=\left( \sum_{j\geq 1} 2^{2(\g - \eta) j}\right)^{1/2}<+\infty$. Note that
the sequence of functions $(\vphi_{j})_{j\geq 1}$ given in (\ref{BaHu}) can be
identified with those given in (\ref{inhpB_norm2}). Hence we may write
$$\sup_{j\geq 1} 2^{\eta j} \|\vphi_{j} *
  f\|_{L^{\infty}}\leq \sup_{j\geq 0} 2^{\eta j} \|\vphi_{j} * f\|_{L^{\infty}},\quad
\vphi_{j}\mbox{ is given by } (\ref{ran_ah}),$$
and then (using (\ref{BaHu})) we obtain:
\begin{equation}\label{sheii}
\|f_{+}\|_{\dot{F}^{\g,a}_{\infty,2}} \leq C
\|f\|_{B^{\eta,a}_{\infty,\infty}}.
\end{equation}
Using (\ref{b_embed}) with $s=2m$, $p=2$, $q=\infty$, $t=\eta$ and
$r=\infty$, we deduce that:
$$B^{2m,a}_{2,\infty}\hookrightarrow B^{\eta,a}_{\infty,\infty}.$$
Therefore, by (\ref{s_embed}), we get
$$W^{2m,m}_{2}\hookrightarrow B^{2m,a}_{2,\infty}\hookrightarrow B^{\eta,a}_{\infty,\infty}$$
which, together with (\ref{sheii}), give:
\begin{equation}\label{result2}
\|f_{+}\|_{\dot{F}^{\g,a}_{\infty,2}} \leq C \|f\|_{W^{2m,m}_{2}}.
\end{equation}
{\bf Step 2} \textit{(An upper bound on
$\|f_{-}\|_{\dot{F}^{-\g,a}_{\infty,2}}$)}. In this step, we will
use the fact that $\partial_{i}g=f_{i}$ (for which we keep denoting
it by $f$, i.e. $f=f_{i}$) for some $i=1,\ldots, n+1$, with $g\in
L^{\infty}(\R^{n+1})$. For $z\in \R^{n+1}$, define
\begin{equation}\label{vhaRmsha}
\Phi(z) = (\partial_{i} \vphi) (z),\quad \vphi \mbox{ is given by
  (\ref{y5tche})},
\end{equation}
and
\begin{equation}\label{vhaR}
\Phi_{j}(z) = 2^{(n+2) j}\Phi(2^{ja}z) \quad \mbox{for all}\quad j\leq
-1.
\end{equation}
Using (\ref{lp:eq2}) we obtain:
\begin{equation}\label{fakahta}
(\partial_{i} \vphi_{j})(z) =
\left\{
\begin{aligned}
& 2^{j} \Phi_{j}(z) \quad &\mbox{if}&\quad i=1,\ldots,n\\
& 2^{2j} \Phi_{j}(z)\quad &\mbox{if}&\quad i=n+1.
\end{aligned}
\right.
\end{equation}
We now compute (see (\ref{conv2}), (\ref{vhaR}) and (\ref{fakahta})):
\begin{eqnarray}\label{Bel2tel}
\|f_{-}\|_{\dot{F}^{-\g,a}_{\infty,2}} &=&  \left\|\left( \sum_{j\leq -1} 2^{-2 \g j}
|\vphi_{j} * f|^{2} \right)^{1/2}\right\|_{L^{\infty}}\\
&\leq& C \sup_{j\leq -1} \|\Phi_{j} * g\|_{L^{\infty}},
\end{eqnarray}
where the constant $C$ is given by:
$$
C^{2}=
\left\{
\begin{aligned}
& \sum_{j\leq -1} 2^{2j(1-\g)} \quad &\mbox{if}&\quad
i=1,\ldots,n\\
& \sum_{j\leq -1} 2^{2j(2-\g)} \quad &\mbox{if}&\quad i=n+1,
\end{aligned}
\right.
$$
which is finite $0<C<+\infty$ under the choice
$$0<\g<1.$$
In order to terminate the proof, it suffices to show that
\begin{equation*}
\|\Phi_{j} * g\|_{L^{\infty}} \leq C \|g\|_{L^{\infty}},
\end{equation*}
which can be deduced, by translation and dilation invariance, from the
following estimate:
\begin{equation}\label{GGLam}
|(\Phi * g)(0)| \leq C \|g\|_{L^{\infty}}.
\end{equation}
Indeed, define the positive radial decreasing function $h(r)=h(\|z\|)$ as follows:
\begin{equation*}\label{sl:eq}
h(r)= \sup_{\|z\| \geq r} |\Phi(z)|.
\end{equation*}
From (\ref{vhaRmsha}), we remark that the function $\Phi$ is the inverse
Fourier transform of a compactly supported function. Hence, we have:
\begin{equation}\label{slom1}
h(0)=\|\Phi\|_{L^{\infty}}<+\infty,
\end{equation}
and the asymptotic behavior
\begin{equation}\label{slom2}
h(r) \leq \frac{C}{r^{n+2}}\quad \mbox{for all}\quad r\geq 1.
\end{equation}
We compute (taking $S_{r}^{n}$ as the $n$-dimensional sphere of radius
$r$):
\begin{eqnarray}\label{770}
|(\Phi * g)(0)| &\leq& \int_{\R^{n+1}} |\Phi(-z)| |g(z)| dz\nonumber\\
&\leq& \int^{\infty}_{0} \left(\int_{S_{r}^{n}}|\Phi(-z)| |g(z)|
  d\sigma(z)\right) dr\nonumber\\
&\leq& C \left(\int^{\infty}_{0} r^{n}h(r) dr\right) \|g\|_{L^{\infty}}.
\end{eqnarray}
Using (\ref{slom1}) and (\ref{slom2}) we deduce that:
\begin{eqnarray*}
\int^{\infty}_{0} r^{n}h(r) dr &=& \int^{1}_{0}r^{n}h(r) dr +
\int_{1}^{\infty}r^{n}h(r) dr\\
&\leq& C\left(\int^{1}_{0}h(0) dr +  \int_{1}^{\infty}
  \frac{r^{n}}{r^{n+2}} dr\right)\\
&\leq& C( \|\Phi\|_{L^{\infty}} + 1)
\end{eqnarray*}
which, together with (\ref{770}), directly implies (\ref{GGLam}). As a
conclusion, we obtain (see (\ref{Bel2tel})):
\begin{equation}\label{result3}
\|f_{-}\|_{\dot{F}^{-\g,a}_{\infty,2}} \leq C \|g\|_{L^{\infty}}.
\end{equation}
{\bf Step 3} \textit{(A lower bound on $\|f\|_{\dot{F}^{0,a}_{\infty,1}}$ and
conclusion)}. Remarking that
$$\|f\|_{L^{\infty}} = \Big\|\sum_{j\in \Z} \vphi_{j} * f
\Big\|_{L^{\infty}} \leq \|f\|_{\dot{F}^{0,a}_{\infty,1}}$$ when
$\hat{f}(0) = 0$, the estimates (\ref{tark}), (\ref{result1}),
(\ref{result2}) and (\ref{result3}) lead directly to the proof.
$\hfill{\Box}$


\section{Proof of Theorem~\ref{theorem-bdd}}\label{sec4} For the sake of simplicity, we only give the
proof in the framework of one spatial dimensions $x=x_{1}$. The
extension to the multi spatial dimensions can be easily deduced and
will be made clear later in this section. Again, the constant $C>0$
that will appear in the following proof may vary
from line to line but will only depend on $m$ and $T$\\

\noindent \textbf{Proof of Theorem~\ref{theorem-bdd}.} We first
remark that the function $f$ can be extended by continuity to the
boundary $\partial \O_{T}$ of $\O_{T}$. Following the same notations
of Ibrahim and Monneau \cite{IM09}, we take $\tilde{f}$ as the
extension of $f$ over
$$\widetilde{\O}_{T} = (-1,2)\times (-T,2T)$$
given by:
\begin{equation}\label{NH}
\tilde{f}(x,t) = \left\{
\begin{aligned}
& \sum_{j=0}^{2m-1} c_{j} f(-\lambda_{j} x, t)\quad
&\mbox{for}&\quad &-1<x<0&,\, &0\leq t\leq T&\\
& \sum_{j=0}^{2m-1} c_{j} f(1+\lambda_{j} (1-x), t)\quad
&\mbox{for}&\quad &1<x<2&,\, &0\leq t\leq T&,
\end{aligned}
\right.
\end{equation}
with $\lambda_{j} = \frac{1}{2^{j}}$ and
$$ \sum_{j=0}^{2m-1} c_{j} (-\lambda_{j})^{k} = 1\quad \mbox{for}\quad k=0,\ldots, 2m-1.$$
For the extension with respect to the time variable $t$, we use the
same extension (\ref{NH}) summing up only to $m-1$. The above
extension (\ref{NH}) has been made in order to have (see for
instance Evans~\cite{Evans}) $\tilde{f} \in
W^{2m,m}_{2}(\widetilde{\O}_{T})$ and
\begin{equation}\label{Berna}
\|\tilde{f}\|_{W^{2m,m}_{2}(\widetilde{\O}_{T})} \leq C
\|f\|_{W^{2m,m}_{2}(\O_{T})}.
\end{equation}
Now let $\mathcal{Z}_{1} \subseteq \mathcal{Z}_{2}$ be two subsets
of $\widetilde{\O}_{T}$ defined by:
$$\mathcal{Z}_{1} = \{(x,t);\, -1/4<x<5/4 \mbox{ and } -T/4<t<5T/4\}$$
and
$$\mathcal{Z}_{2} = \{(x,t);\, -3/4<x<7/4 \mbox{ and } -3T/4<t<7T/4\}.$$
We take the cut-off function $\Psi\in C^{\infty}_{0}(\R^{2})$,
$0\leq\Psi\leq 1$ satisfying:
\begin{equation}\label{5sar}
\Psi(x,t) = \left\{
\begin{aligned}
& 1 \quad \mbox{for}\quad (x,t)\in \mathcal{Z}_{1}\\
& 0 \quad \mbox{for}\quad (x,t)\in \R^{2}\setminus \mathcal{Z}_{2}.
\end{aligned}
\right.
\end{equation}
From (\ref{Berna}), we easily deduce that $\Psi \tilde{f} \in
W^{2m,m}_{2}(\R^{2})$ and
\begin{equation}\label{NHDN}
\|\Psi \tilde{f}\|_{W^{2m,m}_{2}(\R^{2})} \leq C
\|f\|_{W^{2m,m}_{2}(\O_{T})}.
\end{equation}
Hence we can apply the scalar-valued version of inequality
(\ref{I_inq}) (see Remark~\ref{rem_s}) with $i=1$, i.e.
$\partial_{1} g = f$; the new function (for which we give the same
notation) $f = \Psi \tilde{f}\in W^{2m,m}_{2}(\R^{2})$ and $g\in
L^{\infty}(\R^{2})$ given by
$$g(x,t) = \int_{0}^{x} \Psi(y,t)\tilde{f}(y,t) dy.$$
Since $\Psi \tilde{f}$ is of compact support, and (again by the
extension (\ref{NH}))
$\|\tilde{f}\|_{L^{\infty}(\widetilde{\O}_{T})}\leq C
\|f\|_{W^{2m,m}_{2}(\O_{T})}$, we deduce that
\begin{equation}\label{proG}
\|g\|_{L^{\infty}(\R^{2})} \leq C
\|\tilde{f}\|_{L^{\infty}(\widetilde{\O}_T)}\leq C
\|f\|_{W^{2m,m}_{2}(\O_{T})}.
\end{equation}
Collecting the above arguments (namely (\ref{NHDN}) and
(\ref{proG})) together with the fact that (see Ibrahim and Monneau
\cite{IM09})
$$\|\Psi \tilde{f}\|_{BMO^{a}(\R^{2})} \leq C \|f\|_{\overline{BMO}^{a}(\O_T)},$$
inequality (\ref{I_inq}) gives:
\begin{eqnarray*}
\|f\|_{L^{\infty}(\O_{T})} \leq \|\Psi
\tilde{f}\|_{L^{\infty}(\R^{2})} &\leq& C \left(1 + \|\Psi
\tilde{f}\|_{BMO^{a}(\R^{2})}\left(\log^{+}(\|\Psi \tilde{f}\|_{
W_{2}^{2m,m}(\R^{2})} +
\|g\|_{L^{\infty}(\R^{2})})\right)^{1/2}\right)\\
&\leq& C \left(1 +
\|f\|_{\overline{BMO}^{a}(\O_T)}\left(\log^{+}\|f\|_{
W_{2}^{2m,m}(\O_T)}\right)^{1/2}\right).
\end{eqnarray*}
Notice that in the first line of the above inequalities we have used
that $\Psi =1$ in $\O_T$. $\hfill{\Box}$
\begin{rem}
The inequality (\ref{I_inq}) used in the previous proof could have
also been used for $i=2$. In this case we simply take $g(x,t) =
\int_{0}^{t} \Psi(x,s)\tilde{f}(x,s) ds.$
\end{rem}
\begin{rem}
In the case of multi spatial dimensions $x_{i}$, $i=1,\ldots,n$, we
simultaneously apply the extension (\ref{NH}) to each spatial
coordinate while fixing all the other coordinates including time
$t$. However, the extension with respect to $t$ is kept unchanged.
\end{rem}

%
%
%
%

\section{Comparison between parabolic logarithmic inequalities}\label{sec5}
In this section we give the proof of Theorem~\ref{theorem3}.
Throughout all this section, we only consider isotropic function
spaces, i.e. $a=(1,\ldots,1)\in \R^{n+1}$. We only deal with the
parabolic function space $W^{2m,m}_{2}$. As usual, the constant
$C=C(m,n)>0$ may differ from line to line. First of all, we remark
that estimate (\ref{sec5:eq1}) turns out to be true (using the
trivial identity $x^{1/2} \leq 1+x$) if $A:=\|g\|_{L^{\infty}} \leq
C$ for $\|f\|_{W^{2m,m}_{2}} \leq 1$, or if
$B:=\|g\|_{L^{\infty}}/\|f\|_{W^{2m,m}_{2}}\leq C$ for
$\|f\|_{W^{2m,m}_{2}} \geq 1$. This will be proved in the
forthcoming arguments. We start with the following lemmas:
\begin{lemma}\label{sec5:lem1}
Let $n=1,2,3$, $s=\frac{n+1}{2}$, and $m\in \N^{*}$ satisfying
$2m>\frac{n+2}{2}$. For any $g\in L^{2}(\R^{n+1})$ with $f = \nabla
g\in W^{2m,m}_{2}(\R^{n+1})$, we have $g\in \dot{H}^{s}(\R^{n+1})$
and
\begin{equation}\label{sec5:eq2}
\|g\|_{\dot{H}^{s}} \leq \|f\|_{W^{2m,m}_{2}}.
\end{equation}
The norm in the homogeneous Sobolev space $\dot{H}^{s}$ is given by
$\|f\|_{\dot{H}^{s}}^{2}=\int_{\R^{n+1}}
\|\xi\|^{2s}|\hat{f}(\xi)|^{2} d\xi$ where $\|\cdot\|$ is the usual
Euclidean distance.
\end{lemma}
\textbf{Proof.} Follows directly since $1\leq s \leq m$, using the
definition of the norm in $\dot{H}^{s}$. $\hfill{\Box}$
\begin{lemma}\label{sec5:lem2}
Under the same hypothesis of Theorem~\ref{theorem3}, we have:
\begin{equation}\label{sec5:eq3}
\|g\|_{L^{\infty}} \leq C \left[ 1 + \|f\|_{W^{2m,m}_{2}} \left(
\log\left(e + \frac{\|f\|_{W^{2m,m}_{2}} +
\|g\|_{L^{\infty}}}{\|f\|_{W^{2m,m}_{2}}}
\right)\right)^{1/2}\right].
\end{equation}
\end{lemma}
\textbf{Proof.} We consider the isotropic ($a=(1,\ldots,1)$)
homogeneous dyadic partition of unity $(\psi_{j})_{j\in \Z}$ with
$\sum_{j\in \Z} \psi_{j} = 1$ and $\hat{\varphi}_{j} = \psi_{j}$.
Fix some $0<\g<1$, and take an arbitrary $N\in \N^*$. We write:
\begin{equation}\label{sec5:eq4}
\|g\|_{L^{\infty}} \leq \sum_{j\leq 0}\|\vphi_{j}*g\|_{L^{\infty}} +
\sum_{j=1}^{N} \|\vphi_{j} * g\|_{L^{\infty}} + \sum_{j>N}
\|\vphi_{j} * g\|_{L^{\infty}}.
\end{equation}
We estimate the right-hand side of (\ref{sec5:eq4}). Benstein's
inequality gives:
\begin{equation}\label{sec5:eq5}
\|\vphi_{j} * g\|_{L^{\infty}} \leq C
2^{\left(\frac{n+1}{2}\right)j} \|\vphi_{j} * g\|_{L^2}.
\end{equation}
We let $s=\frac{n+1}{2}$. Using (\ref{sec5:eq5}), we compute:
\begin{eqnarray}\label{sec5:eq0}
\sum_{j\leq 0}\|\vphi_{j}*g\|_{L^{\infty}} &\leq& C\sum_{j\leq 0}
2^{sj} \|\vphi_{j} * g\|_{L^2}\nonumber\\
&\leq& C\sum_{j\leq 0}
2^{sj} \|\widehat{\vphi_{j} * g}\|_{L^2}\nonumber\\
&\leq& C\sum_{j\leq 0}
2^{sj} \|\psi_{j} \hat{g}\|_{L^2}\nonumber\\
&\leq& \frac{C}{1 - 2^{-s}} \|g\|_{L^{2}}.
\end{eqnarray}
Again, using (\ref{sec5:eq5}), we obtain:
\begin{eqnarray*}
\sum_{j=1}^{N} \|\vphi_{j} * g\|_{L^{\infty}} &\leq& C
\sum_{j=1}^{N}
2^{sj} \|\vphi_{j} * g\|_{L^{2}}\\
&\leq& C N^{1/2} \left( \sum_{j=1}^{N} 2^{2sj} \|\vphi_{j} *
g\|^{2}_{L^{2}}\right)^{1/2}\\
&\leq& C N^{1/2} \|g\|_{\dot{B}^{s}_{2,2}},
\end{eqnarray*}
which, together with the fact that $\dot{B}^{s}_{2,2} \simeq
\dot{H}^{s}$, and estimate (\ref{sec5:eq2}) of
Lemma~\ref{sec5:lem1}, yield:
\begin{equation}\label{sec5:eq6}
\sum_{j=1}^{N} \|\vphi_{j} * g\|_{L^{\infty}} \leq  C N^{1/2}
\|f\|_{W^{2m,m}_{2}}.
\end{equation}
The last term of the right-hand side of (\ref{sec5:eq4}) can be
estimated as follows:
\begin{eqnarray}\label{sec5:snh}
\sum_{j>N} \|\vphi_{j} * g\|_{L^{\infty}} &=& \sum_{j>N}2^{-j\g} (
2^{j\g} \|\vphi_{j}* g\|_{L^{\infty}})\nonumber\\
&\leq& \big(\sum_{j>N}2^{-j\g}\big) \sup_{j\in \Z} 2^{j\g}
\|\vphi_{j}* g\|_{L^{\infty}}\nonumber\\
&\leq& 2^{-\g N} \left(\frac{2^{-\g}}{1- 2^{-\g}}\right)
\|g\|_{\dot{B}^{\g}_{\infty,\infty}}.
\end{eqnarray}
We know that $\dot{B}^{\g}_{\infty,\infty} \simeq \dot{C}^{\g}$; the
homogeneous H\"{o}lder space whose semi-norm can be estimated as
follows:
$$\|g\|_{\dot{C}^{\g}} = \sup_{z_{1}\neq z_{2}}
\frac{|g(z_{1}) - g(z_{2})|}{\|z_{1} - z_{2}\|^{\g}} \leq
\|f\|_{W^{2m,m}_{2}} + \|g\|_{L^{\infty}}.$$ This, together with
(\ref{sec5:snh}) yield:
\begin{equation}\label{sec5:eq7}
\sum_{j>N} \|\vphi_{j} * g\|_{L^{\infty}} \leq C 2^{-\g N}
(\|f\|_{W^{2m,m}_{2}} + \|g\|_{L^{\infty}}).
\end{equation}
Combining (\ref{sec5:eq4}), (\ref{sec5:eq0}), (\ref{sec5:eq6}) and
(\ref{sec5:eq7}), we finally get:
$$\|g\|_{L^{\infty}} \leq C \big(1 + N^{1/2} \|f\|_{W^{2m,m}_{2}} + 2^{-\g N} (\|f\|_{W^{2m,m}_{2}} +
\|g\|_{L^{\infty}})\big).$$ By optimizing (as in Step 2 of
Lemma~\ref{lem2}) in $N$ the above inequality, the proof easily
follows.$\hfill{\Box}$\\

We are now ready to give the proof of Theorem~\ref{theorem3}.\\

\noindent \textbf{Proof of Theorem~\ref{theorem3}.} As it was
already mentioned in the beginning of this section, the proof
relies on considering two cases.\\

\noindent \textbf{Case 1} ($\|f\|_{W^{2m,m}_{2}}\leq 1$). Let
$A:=\|g\|_{L^{\infty}}$. Using inequalities (\ref{nbbo2}) and
(\ref{sec5:eq3}), we obtain:
$$A \leq C [1 + (\log (e + 1 + A))^{1/2}],$$
which directly implies that:
$$A\leq C,$$
and hence (\ref{sec5:eq1}) is obtained.\\

\noindent \textbf{Case 2} ($\|f\|_{W^{2m,m}_{2}}\geq 1$). Dividing
inequality (\ref{sec5:eq3}) by $\|f\|_{W^{2m,m}_{2}}$, we obtain:
$$\frac{\|g\|_{L^{\infty}}}{\|f\|_{W^{2m,m}_{2}}} \leq C \left[ 1 +  \left(
\log\left(e + 1 + \frac{ \|g\|_{L^{\infty}}}{\|f\|_{W^{2m,m}_{2}}}
\right)\right)^{1/2}\right].$$ Letting
$B:=\|g\|_{L^{\infty}}/\|f\|_{W^{2m,m}_{2}}$, we can easily see that
$B$ satisfies (as the term $A$ in Case~1):
$$B \leq C [1 + (\log (e + 1 + B))^{1/2}],$$
which shows that:
$$B \leq C,$$
and the proof is done. $\hfill{\Box}$

\end{document}